# Geometry of right-angled Coxeter groups on the Croke–Kleiner spaces

**Yulan Qing**[1]

**Abstract** In this paper we study the right-angled Coxeter groups that acts geometrically on the Salvetti complex of a certain right-angled Artin group, which we refer to as Croke–Kleiner spaces. We prove that any right-angled Coxeter group that acts geometrically on the Croke–Kleiner spaces acts with $\pi/2$ angles between reflecting axes, while the quasi-isometric right-angled Artin group can act with angles that are any real number in the range $(0, \pi/2]$. The contrast between the two examples shows that in this case a right-angled Coxeter group is geometrically more "rigid" than its quasi-isometric counterpart.

**Mathematics Subject Classification** 20 · 51 · 57

## 1 Introduction

The motivation of this study is to better understand the relationship between right-angled Coxeter groups and right-angled Artin groups. It is known that every right-angled Artin group embeds as a finite index subgroup of a right-angled Coxeter group [3], so we ask what are the geometric elements, if any, that distinguish a pair of quasi-isometric groups, one from each category. We would like to generalize the result to all pairs of quasi-isometric right-angled Artin groups and right-angled Coxeter groups.

In this paper we provide a pair of quasi-isometric groups: the right-angled Artin group acts on the Croke–Kleiner space [2] geometrically with flexible angles, while the right-angled Coxeter group can only act with fixed angles. Furthermore, we prove that *any* right-angled Coxeter group that acts geometrically on this family of spaces has the same constraint.

✉ Yulan Qing
yulan.qing@gmail.com; yulan.qing@utoronto.ca

[1] Department of Mathematics, University of Toronto, Toronto, ON M5S 2E4, Canada

Therefore, the "right-angled" in the terminology of right-angled Coxeter groups takes on a literal meaning.

The main result is Theorem 1.1. We show that:

**Theorem 1.1** *Let W be a right-angled Coxeter group that acts geometrically on the Croke–Kleiner spaces, then the intersection angle in $T_2$ has to be $\pi/2$.*

The outline of the paper is as follows. In Sect. 2 come preliminaries needed in this paper. In Sect. 3 we introduce the Croke–Kleiner spaces; in Sect. 4 we prove the claim.

## 2 Preliminaries

In this section we give basic definitions and facts concerning CAT(0) geometry and boundaries, all of whose proofs can be found in [1]. We also give the definitions and facts we need concerning right-angled Artin and Coxeter groups.

### 2.1 CAT(0) spaces and their boundaries

A metric space $X$ is $CAT(0)$ if geodesic triangles in $X$ are at least as thin as a triangle in Euclidean space with the same side lengths. It follows immediately from the definition that CAT(0) spaces are uniquely geodesic and thus contractible via geodesic retraction to a base point in the space.

Recall that a metric space $X$ is *proper* if closed metric balls are compact. In this case, $X$ can be compactified via the *visual boundary* of $X$. The points of this boundary are equivalence classes of geodesic rays defined as follows:

A *geodesic ray* in $X$ is a geodesic $c : [0, \infty) \to X$. Consider the set of geodesic rays in $X$. Two geodesic rays $c_1$ and $c_2$ are said to be *asymptotic* if $f(t) := d(c_1(t), c_2(t))$ is a bounded function. The set of equivalence classes is denoted by $\partial X$. If $\xi \in \partial X$ and $c$ is a geodesic ray belonging to $\xi$, we write $c(\infty) = \xi$.

The following is a lemma in CAT(0) geometry:

**Lemma 1** *For any $\xi \in \partial X$ and any $x \in X$, there is a unique geodesic ray $c_{x\xi} : [0, \infty) \to X$ with $c_{x\xi}(0) = x$ and $c_{x\xi}(\infty) = \xi$. The image of $c_{x\xi}$ is denoted by $x\xi$.*

The topology on $\partial X$ called *cone topology* has as a basis the open sets of $X$ together with the sets

$$U(x, \xi, R, \epsilon) = \{\xi' \in \partial X | d(c_{x\xi'}(R), c_{x\xi}(R)) < \epsilon\}$$

where $x \in X$, $\xi \in \partial X$ and $R > 0$, $\epsilon > 0$. The topology on $X$ induced by the cone topology coincides with the metric topology on $X$.

This topology looks as if it depends on the base-point $x$ in the above description of open sets, however the previous lemma shows that there is a natural change of base-point homeomorphism when the base-point is changed (Fig. 1).

The set $\partial X$ together with the cone topology is called the *visual boundary* of $X$, denoted $\partial_\infty X$. It is also referred to as *ideal boundary*.

### 2.2 Quasi-isometry and quasi-isometric embeddings

**Definition 2.1** Let $(X_1, d_1)$ and $(X_2, d_2)$ be metric spaces. A (not necessarily continuous) map $f : X_1 \to X_2$ is called a $(\lambda, \epsilon)$-*quasi-isometric embedding* if there exist constants

**Fig. 1** A basis for open sets

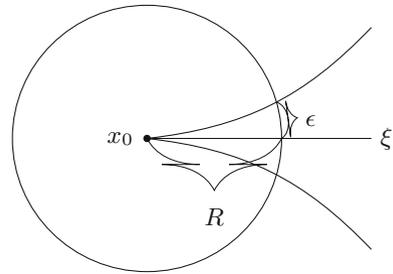

$\lambda \geq 1$ and $\epsilon \geq 0$ such that for all $x, y \in X_1$

$$\frac{1}{\lambda}d_1(x, y) - \epsilon \leq d_2(f(x), f(y)) \leq \lambda d_1(x, y) + \epsilon$$

If, in addition, there exists a constant $C \geq 0$ such that every point of $X_2$ lies in the $C$-neighborhood of the image of $f$, then $f$ is called a $(\lambda, \epsilon)$-quasi-isometry. When such a map exists, $X_1$ and $X_2$ are said to be *quasi-isometric*.

### 2.3 Right-angled groups

**Definition 2.2** Let $\Gamma$ be a finite, simplicial graph. The vertex set is denoted by $V$ and edge set $E$ denotes a set of unordered pairs of vertices. The requirement of being simplicial means the diagonal $V \times V$ is excluded from the edge set. The *right-angled Artin group* on $\Gamma$ is the group

$$A(\Gamma) \cong \rangle V | [v_i, v_j] = 1 \text{whenever} (v_i, v_j) \in E \langle$$

$A(\Gamma)$ is generated by the vertices of $\Gamma$, and the only relations are given by commutation of adjacent vertices.

**Definition 2.3** A *right-angled Coxeter Group* on a finite, simplicial $Gamma = (V, E)$ has the presentation:

$$A(\Gamma) \cong \rangle V | v_i^2 = 1 \forall i, [v_i, v_j] = 1 \text{whenever} (v_i, v_j) \in E \langle$$

Just as for right-angled Artin groups, the presentation for a right-angled Coxeter group can be given by a finite simplicial graph with the understanding that each vertex now represents a generator of order 2.

Here are more basic facts of right-angled Coxeter groups [5]:

- If $s_i$ is not adjacent to $s_j$, then the order of $s_i s_j$ is infinite.
- A right-angled Coxeter group is abelian if and only if it is finite which is true if and only if the defining graph is complete.
- A right-angled Coxeter group $W$ has a non-trivial center if and only if it can be written as $W' \times \mathbb{Z}_2$ for a right-angled Coxeter group $W'$.

## 3 Croke–Kleiner spaces and their boundaries

It is well-known that if $A$ is a right-angled Artin group, then $A$ acts geometrically on a CAT(0) cube complex, namely the universal cover of its Salvetti complex. We will construct

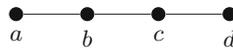

**Fig. 2** The defining graph

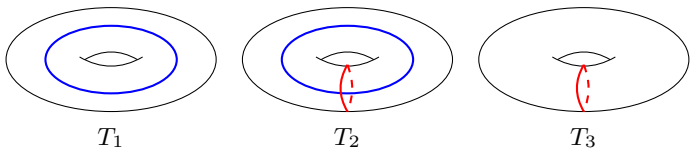

**Fig. 3** Tori complex

the universal cover of the Salvetti complex of a specific right-angled Artin group in this section (Fig. 2). The group A has as its defining graph:

Here is the construction laid out by Croke and Kleiner [2]. Start with a flat torus $T_2$ with the property that a pair $b$, $c$ of unoriented, $\pi_1$-generating simple closed curves in $T_2$ meets at a single point at an angle $\theta_2 = \frac{\pi}{2}$. Let $b, c$ both have length 1. Let $T_1, T_3$ be flat tori containing simple closed essential loops, $a, b_1$ and $c_1, d$, respectively, such that $length(b_1) = length(b), length(c_1) = length(c)$. Let $a, d$ also have length 1. Likewise $\theta_1$ and $\theta_3$ denote the intersecting angles of the generating loops in $T_1$ and $T_3$, respectively. Let $Y$ be the union of $T_1, T_2, T_3$ with $b_1$ identified isometrically with $b$ and $c_1$ with $c$. Let $X$ be the universal cover of $Y$. Let $Y_1 = T_1 \cup T_2$, and let $Y_2 = T_3 \cup T_2$. That $X$ is CAT(0) cube complex follows from the Equivariant Gluing Theorem 11.18 [1] (Fig. 3).

We obtain an uncountably infinite family of CAT(0) spaces by changing the geometry of $X$ in such a way that it is no longer cubical yet it is still CAT(0) and the group A still acts on the new spaces geometrically. Specifically, we can change the angle $\theta_2$ to be any real number $0 < \theta_2 \leq \frac{\pi}{2}$. This particular change of geometry was studied in the original paper [2] where they proved that changing the angle from $\frac{\pi}{2}$ to any other value changes the homeomorphism type of $\partial_\infty X$. This was further investigated by Wilson [9] where it was shown that any two different angles give different visual boundaries.

We can also change the lengths of $a, b, c, d$, which affects the translation lengths of the actions. The *geometric data* associated with a Croke–Kleiner space consists of three intersecting angles and four translation distances. The three intersection angles are that of the intersecting angle of the three pairs of $\pi_1$-generating, simple closed curves on the three tori, which we denote $\theta_1, \theta_2, \theta_3$. The four lengths are the translation distances of $a, b, c, d$. All spaces constructed this way, with different geometric data, are quasi-isometric and each will be referred to as a *Croke–Kleiner space*.

We now describe the structure of any Croke–Kleiner space $X$.

**Definition 3.1** A *barrier* in $X$ is a maximal connected component of the universal cover of $T_2$. A *block* in $X$ is a maximal, connected component of the universal cover of $Y_i$, denoted $X_i$.

Each block, as well as each barrier, is a closed, connected and locally convex subset of $X$. Let $\mathscr{B}$ denote the collection of all blocks and $\mathscr{W}$ the collection of all barriers. We prove later that $\mathscr{B}$ and $\mathscr{W}$ are countably infinite sets.

Let $\mathscr{T}_4$ be the regular, 4-valent, infinite tree that is graph isomorphic to the Cayley graph of $F_2$ with two generators. A block is isometric to the metric product of $\mathscr{T}_4$ with appropriate edge lengths with the real line $\mathbb{R}$. The intersection of two blocks can be either an empty set or a barrier. Two blocks are *adjacent* if their intersection is a barrier (Fig. 4).

**Fig. 4** A block

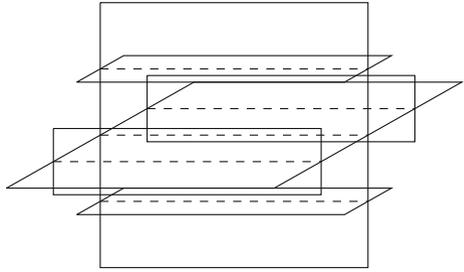

**Fig. 5** $\mathcal{T}_t$

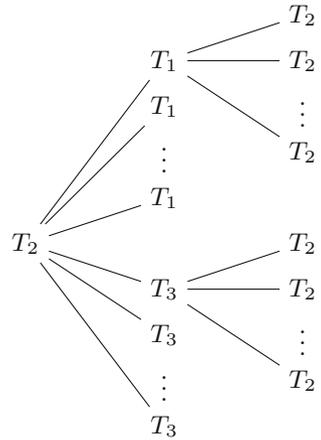

There are two trees on which a Croke–Kleiner space can be projected. $\mathcal{T}_b$ denotes the tree whose vertices are blocks and whose edges denote the adjacency of pairs of blocks. $\mathcal{T}_t$ denotes the tree whose vertices are lifts of $T_1$, $T_2$ and $T_3$, and two vertices are adjacent if and only if two respective planes have nonempty intersection. Both trees are locally countably infinite, infinite trees, with different labeling. The figure shows $\mathcal{T}_t$. Here each vertex is labeled by the torus they are a lift of. $T_2$'s are adjacent to countably infinite $T_1$'s and $T_3$'s; while a $T_1$(or $T_3$) only has $T_2$'s as its neighboring vertices (Fig. 5).

$\mathcal{T}_t$ and $\mathcal{T}_b$ corresponds to different Bass–Serre decomposition of the group

$$\langle a, b, c, d | [a, b], [b, c], [c, d] \rangle.$$

They are the *Bass–Serre trees*.

## 4 Proof of the main theorem

In this section we prove the main result of this paper, Theorem 1.1.

In the rest of the paper , let $W$ be any right-angled Coxeter group with the standard presentation:

$$W = \{s_1, s_2, \ldots s_n | s_i^2 \text{ for all } i, [s_i, s_j] \text{ for some pairs } i, j\}$$

The main theorem of this paper is:

**Theorem 4.1** *Suppose a right-angled Coxeter group W acts geometrically on a Croke–Kleiner space, then the intersection angle on the middle torus must be $\frac{\pi}{2}$.*

Recall a *special flat* is a flat that is a lift of $T_i$. We start with a lemma about the action of the generators and the stabilizers on each special flat.

**Lemma 4.2** *W preserves the set of special flats.*

*Proof* Since the actions of the group elements are isometric, they induce homeomorphisms on the visual boundary. Homeomorphisms of a visual boundary preserve block boundaries [2]. Hence the homeomorphisms also preserve the intersections of block boundaries, which are circles. Circles are boundaries of special flats. Each circle corresponds to two points in the tree factor of a block, and since the trees are hyperbolic, two points in the tree factor picks out exactly one plane that is in the equivalence class of the special flat determined group-theoretically by the two blocks, which has to be special flat itself. Therefore the actions preserve the special flats. □

**Corollary 4.3** *For Croke–Kleiner complexes, each generator of W acts on the nerve tree without inversion.*

*Proof* By construction, each edge in the nerve is incident to a pair of vertices, one of them denotes a plane that is a special flat, the other denotes a plane that is a lift of either $T_1$ or $T_3$. Since the actions preserve special flat, they cannot invert any edges. □

Since the action is without inversion, we use a Bass–Serre Theorem [8, p. 39]

**Theorem 4.4** *If W is a group acting without inversion on tree X, then there is a finite tree T whose vertex groups are stabilizers of the vertices and whose edge groups are stabilizers of the edges and the group W is the Bass–Serre group of that construction.*

How do the generators of W act on the finite tree T in the previous theorem? By construction the finite tree T is the strict fundamental domain of W. By the covering property of fundamental domains, the fixed point set of $s_i$ must intersect a translate of T. If the intersection is more than one point, it violates the definition of a strict fundamental domain, then thus intersection must be a point. But since the action is without inversion, then the fixed point set must be a vertex, not a midpoint of an edge. Therefore each generator of W (or its appropriate translate) stabilizes a vertex in T.

In general, if a group element stabilizes a vertex in T, we would like to show that it is a conjugate of a word in the subset of generators that stabilize that vertex in T. Here is the gist of the argument. Go to a subsequence that does not fix the vertex and examine the image of T under the actions of each generator in that subsequence. Then the image of that vertex forms a "circle" which contradicts the fact it is embedded in a tree(the nerve tree), which means that it must be reflected "out and back" at some point, i.e. some of the images of T overlap, which again contradicts the definition of a strict fundamental domain. Therefore the only possibility is that the images of that vertex under each generators in the subsequence is itself, which is our claim.

**Lemma 4.5** *Suppose a right-angled Coxeter group W acts geometrically on the Croke–Kleiner space. If a group element w fixes a special flat $\widetilde{T}_i$ set-wise, suppose $w = s_k s_{k-1} \ldots s_2 s_1$, then each $s_i$ fixes $\widetilde{T}_i$ set-wise.*

*Proof* Without loss of generality, suppose $s_1$ does not fix the special flat $T_i$, otherwise let $w = s_k s_{k-1} \ldots s_2$. Let $j$ be the smallest number such that the sub-word $s_j s_{j-1} \ldots s_2 s_1$ fixes the $T_i$. Consider generators $s_1$ and $s_j$. In $\mathcal{T}_0$, $s_j$ and $s_1$ each label a vertex, $v_{s_j}$ and $v_{s_1}$. $\mathcal{T}_0$ also contains a lift of $T_i$, label it $v_0$. Since $\mathcal{T}_0$ is a tree, there are unique paths $(v_0, v_{s_j})$ and $(v_0, v_{s_1})$. The word $s_{j-1} s_{j-2} \ldots s_2$ takes the edges $(v_0, v_{s_j})$ to the edges $(v_0, v_{s_1})$. This contradicts the assumption that $T_2$ is a strict fundamental domain. Therefore, each $s_i$ fixes $\widetilde{T_i}$ set-wise. □

Next we prove that stabilizer subgroups act geometrically.

**Proposition 4.6** *Given the universal cover of $T_i$, denoted $\widetilde{T_i}$ consider the stabilizer subgroup $Stab(\widetilde{T_i})$, then $Stab(\widetilde{T_i})$ is generated by a (conjugate) of a subset of the generating set $\{s_1, s_2, \ldots s_n\}$, respectively.*

*Proof* Each generator acts simplicially on the nerve tree of special flats. Furthermore, let every edge has length 1, then each group element acts isometrically on the tree. Each generator is of order two. Therefore the fixed point set of each generator acting on the nerve tree is either an induced subgraph or the midpoint of an edge. Lemma 4.3 rules out the latter case. By Corollary 4.4, there exists a minimal finite tree that is the strict fundamental domain of $W$ on the nerve, which we denote by $\mathcal{T}$. This tree is the fundamental domain of the group acting on this tree, therefore generators of $W$ does not take points of $\mathcal{T}$ to points of $\mathcal{T}$. By Lemma 4.5, each group element that stabilizes a vertex of $\mathcal{T}_0$ is generated by a subset of generators that stabilizes the vertex. Thus $Stab(\widetilde{T_i})$ are *special subgroups*, i.e. they are generated by a subset of generators. □

We claim the group $Stab(\widetilde{T_i})$ acts on $\widetilde{T_i}$ properly discontinuously, cocompactly and by isometries.

**Proposition 4.7** *Given that $W$ acts geometrically on the CK space, $Stab(\widetilde{T_i})$ acts properly discontinuously, cocompactly and by isometries on the special flat.*

*Proof* The group acts by isometries because it acts by isometries on the whole space. the group acts cocompactly on the the space for the following reason. If $K$ is a fundamental domain for $W \curvearrowright X$, then

$$K \cap \widetilde{T_i}$$

is the fundamental domain for the actions of $Stab(\widetilde{T_i})$. Therefore $Stab(\widetilde{T_i})$ acts cocompactly on the special flat it stabilizes.
The group acts properly discontinuously because by Bass–Serre theory, the stabilizer subgroups contains all the relators that has to do with the elements—part of the consequence of being a special subgroup. So the properly discontinuously carries down, also undistorted. □

Next we study a right-angled Coxeter group acting cocompactly and by isometries on a 2-dimensional Euclidean special flat. For a presentation of a right-angled Coxeter group

$$W = \{s_1, s_2, \ldots s_n | s_i^2 \text{ for all } i, [s_i, s_j] \text{ for some pairs } i, j\}$$

Consider the defining graph of the group. First one can rule out the defining graphs on less than or equal to three vertices since they either have 0, 2, or infinitely many ends. Indeed, the number of ends of a group is a quasi-isometry invariant and the plane has one end so a group with 0 or more than one end cannot act geometrically on the plane by the Svarc–Milnor Lemma.

Recall Gromov's Theorem [1,4]:

**Theorem 4.8** *If a finitely generated group is quasi-isometric to $\mathbb{Z}^n$ then it contains $\mathbb{Z}^n$ as a subgroup of finite index.*

**Lemma 4.9** *(Key Lemma) Suppose W is a right-angled Coxeter group acting cocompactly and by isometries on the special flat $\mathbb{E}^2$. Then we claim that W must be the direct product of two copies of the infinite dihedral group.*

*Proof* We know that the group $W$ has at least four generators. Since $W$ contains $\mathbb{Z}^2$ as a subgroup of finite index, it is not hyperbolic. By [7], if $\Gamma$ is the defining graph of $W$, then in $\Gamma$ there exists induced subgraphs $A$, $B$ such that $\langle A \rangle$, $\langle B \rangle$ are infinite and $A *_{join} B$ is a subgraph of $\Gamma$. In particular, there exists two infinite order elements $\gamma_1' = s_1 t_1$, $\gamma_2' = s_2 t_2$ such that the subgraph on the vertices $s_1, s_2, t_1, t_2$ is a join of two pairs of non-adjacent vertices. The subgraph is a chord-less 4-cycle, where $s_1$ is adjacent to $s_2$ and $t_2$, and $t_1$ is adjacent to $s_2$ and $t_2$.

The actions of $s_1, t_1, s_2, t_2$ are order-2 isometries of the special flat, which are either reflecting across a straight line $l$, or rotating around a point $p$ by $\pi$. Two such elements commute in the following cases:

- $l_1$ and $l_2$ intersecting at right angle
- $l \cap p \neq \phi$

An infinite order action must be a composition of these order-2 isometries as one of the following cases:

1. $l \cap p = \phi$
2. $l_1 \cap l_2 = \phi$
3. $p_1 \cap p_2 = \phi$

In the $\mathbb{Z}^2$ subgroup, there are two elements of infinite order, both generators in one of the three pairs of elements commutes with both generators of another one, not necessarily different, of the three pairs.

(1) and (1), impossible since a point cannot simultaneously coincide with another point off the line and be on the line, For the same reason, (1) and (3) is also impossible.

(1) and (2), impossible since there is only one straight line that passes perpendicularly through another line and a point off that line;

(2) and (3), impossible, since one point cannot be on two parallel lines;

(3) and (3), impossible since one point cannot coincides with two points.

Therefore the only possibility is (2) and (2): two pairs of parallel lines intersecting at right angle. The defining of this group consists of four vertices and four edges connecting up to a four-gon.

To have this group as a subgroup of finite index, by the Finite Index Lemma [6] we must have in the defining graph a complete graph joined to the chord-less 4-cycle. This is to say the generators not in the chord-less 4-cycle commutes with the four reflections. By the previous argument, there cannot be order-2 symmetries of the special flat that commutes with all four reflections. Therefore, if a right-angled Coxeter group acts geometrically on a special flat, the actions of the group restricted to the special flat is isomorphic to

$$G = D_\infty \times D_\infty = \langle a, b, c, d | a^2, b^2, c^2, d^2, [a,c], [a,d], [b,c], [b,d] \rangle$$

□

A direct corollary is the following,

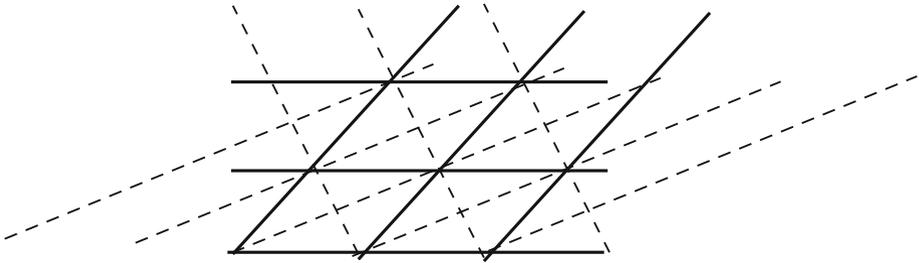

**Fig. 6** Gluing theorem

**Corollary 4.10** *If a right-angled Coxeter group acts geometrically on a special flat, then it is isomorphic to $D_\infty \times D_\infty$. The $D_\infty \times D_\infty$ acts on the special flat like two pairs reflections cross lines. Each pair consists of two reflections whose fixed-point sets are parallel axes, and the two pairs of axes intersect at right angle.*

Next we study how the stabilizer subgroups piece together and determine the gluing angle of the complex.

**Theorem 4.11** *If a right-angled Coxeter group acts geometrically on the Croke–Kleiner space and preserves special flats, then the angle $\theta_2$ must be $\pi/2$.*

*Proof* Consider a special flat of type $T_1$ or $T_3$, without loss of generality, let it be $T_1$. Each of these flats is adjacent to countably many lifts of $T_2$. The intersections of a flat of type $T_i$ with a flat of type $T_j$ is labeled $l_{i,j}$. Any action on $T_1$ preserves the set of all $l_{2,1}$, therefore the axes of the generators are parallel and perpendicular to the $l_{2,1}$. Moreover, two of the generators reflect $l_{2,1}$ across a point.

Now consider the special flats of type $T_2$. In these flats, there are two sets of intersections with neighboring special flats, labeled accordingly $l_{2,1}$ and $l_{2,3}$. All the $l_{2,1}$s are parallel to one another; all the $l_{2,3}$s are parallel to one another. Consider the angle $\theta$ between $l_{2,1}$ and $l_{2,3}$. Suppose $\theta \neq \pi/2$, then the only possibility for a set of four reflection, configured in the way specified in Key Lemma, can take intersections to intersections is to have them reflect across the diagonals of the unit parallelograms in the special flat, as shown in Fig. 6. In Fig. 6, the solid lines are $l_{2,1}$ and $l_{2,3}$, the dashed lines are the axes of reflections.

In this case, it takes a two-letter word to reflect $l_{2,1}$ onto itself across a point. We argued in the first paragraph that there are generators that reflect $l_{2,1}$ to itself across a point. Since $l_{2,1}$ is also in the flat that is a lift of $T_2$, the same generator then act as reflection on the corresponding $T_2$ and its axis intersects $l_{2,1}$. However there are already reflection axes intersecting $l_{2,1}$ as established in the previous paragraph and neither of them reflect $l_{2,1}$ onto itself. Therefore we need to have a third reflection axis that is not parallel to the two existing axes. This configuration contradicts the Key Lemma. Therefore, it is not possible to have the intersection angle of $l_{2,1}$ and $l_{2,3}$ be $\theta \neq \pi/2$. $\square$